\documentclass[twoside,reqno]{amsart}

\usepackage{amsmath}
\usepackage{amssymb,latexsym,amstext,amsthm}
\usepackage{verbatim} 
\usepackage{eucal} 

\theoremstyle{plain}
\newtheorem{thm}[equation]{Theorem}
\newtheorem{pro}[equation]{Proposition}
\newtheorem{cor}[equation]{Corollary}
\newtheorem{lem}[equation]{Lemma}

\theoremstyle{definition}

\newtheorem{con}[equation]{Convention}
\newtheorem{DEF}[equation]{Definition}
\newtheorem{rem}[equation]{Remark}

\def\0b{\bar{0}}

\def\epb{\bar{\ep}}

\def\hz{\hat{z}}

\def\andd{\quad\hbox{and}\quad}
\def\ind{\hbox{ind}}

\def\v{{\mathcal V}}

\def\vd{\dot{\mathcal V}}
\def\vz{{\mathcal V}^{0}}
\def\vt{\tilde{\mathcal V}}

\def\fm{(\cdot,\cdot)}
\def\a{\alpha}

\def\ac{\alpha^{\vee}}
\def\bc{\betta^{\vee}}

\def\w{{\mathcal W}}

\def\sub{\subseteq}
\def\rd{\dot{R}}

\def\lam{\lambda}

\def\1k{\frac{1}{k}}

\def\la{\langle}
\def\ra{\rangle}
\def\rds{\dot{R}_{sh}}
\def\rdl{\dot{R}_{lg}}

\def\rs{R_{sh}}
\def\rl{R_{lg}}

\def\d{\delta}

\def\b{\beta}
\def\bc{\beta^{\vee}}

\def\sg{\sigma}

\def\quadd{\quad\quad}

\def\bbbz{{\mathbb Z}}
\def\bbbr{{\mathbb R}}

\def\B{\mathcal B}

\def\ep{\epsilon}

\def\span{\hbox{span}}

\def\proof{\noindent{\bf Proof. }}

\def\bi{{\bf I}}
\def\bii{{\bf II}}

\def\hw{\hat{w}}
\def\hz{\hat{z}}

\def\hu{\hat{u}}

\def\th{\hat{t}}

\def\a{\alpha}

\def\ac{\alpha^{\vee}}

\def\andd{\quad\hbox{and}\quad}

\def\b{\beta}

\def\bc{\beta^{\vee}}

\def\d{\delta}

\def\bbbr{{\mathbb R}}

\def\ep{\epsilon}
\def\ind{\hbox{ind}}
\def\fm{(\cdot,\cdot)}

\def\sub{\subseteq}
\def\rd{\dot{R}}

\def\lam{\lambda}

\def\1k{\frac{1}{k}}

\def\la{\langle}
\def\ra{\rangle}
\def\rds{\dot{R}_{sh}}
\def\rdl{\dot{R}_{lg}}

\def\rs{R_{sh}}
\def\rl{R_{lg}}

\def\GL{GL}

\def\qed{\hfill$\Box$\vspace{5mm}}
\def\sg{\sigma}

\def\rtimes{R^{\times}}

\def\quadd{\quad\quad}

\def\vd{\dot{\mathcal V}}

\def\vz{{\mathcal V}^{0}}

\def\vt{\tilde{\mathcal V}}

\def\w{{\mathcal W}}

\def\bbbz{{\mathbb Z}}

\def\r{R}

\def\supp{\hbox{supp}}

\def\esupp{\hbox{Esupp}}

\def\tJ{\tau_{_J}}

\def\J{\mathcal J}

\def\matB{\mathcal B}

\begin{document}
\setcounter{page}{1} \setcounter{page}{1}
\title{extended affine Weyl groups: \\Presentation by conjugation via\\ integral collection}

\author{Saeid Azam, Valiollah Shahsanaei}
\address
{Department of Mathematics\\ University of Isfahan\\Isfahan, Iran,
P.O.Box 81745-163} \email{azam@sci.ui.ac.ir, saeidazam@yahoo.com.}
\address{Department of Mathematics\\ University of Qom\\Qom, Iran, P.O.Box 3716146611}\email{shahsanaei@qom.ac.ir.}
\thanks{The authors would like to thank the Center of Excellence for Mathematics, University of Isfahan.}

\begin{abstract}
We give several necessary and sufficient conditions for the existence of {\it the presentation by conjugation} for a non-simply laced extended affine Weyl group. We invent a computational tool by which one can determine simply the existence of the presentation by conjugation for an extended affine Weyl group. As an application, we determine the existence of the presentation by conjugation for a large class of extended affine Weyl groups.
\end{abstract}
\maketitle

\markboth{PRESENTATION BY CONJUGATION}{S. AZAM, V. SHAHSANAEI}

\medskip
\setcounter{section}{-1}
\section{\bf Introduction}\label{introduction}
Let $R$ be an
extended affine root system, then its Weyl group $\w$ is said to have the {\it presentation by conjugation} if it has the following
presentation:

generators: $\hat{w}_\a$, $\a\in\rtimes$,

relations: $\hat{w}_\a^2=1,\quad$ $\hat{w}_\a\hat{w}_\b\hat{w}_\a=
\hat{w}_{w_\a(\b)}$,$\quad (\a,\b\in\rtimes),$

\noindent where $\rtimes$ is the set of non-isotropic roots of $R$ and $w_\a$ is the reflection based on $\a\in\rtimes$. In other
words $\w$ has the presentation by conjugation, if it is isomorphic to the presented group $\hat{\w}$, defined by the above
generators and relations. Since $\w$ is a Hopfian group (see Lemma \ref{AA}), it follows that $\w$ has the presentation by conjugation if and only if the epimorphism $\psi:\hat{\w}\rightarrow\w$ induced from the assignment $\hw_\a\mapsto w_\a$, $\a\in\rtimes$ is one-to-one.

Let us recall briefly what is known about this presentation.
The following subclasses of extended affine Weyl groups are known to have the presentation by conjugation:

\begin{itemize}
\item{} finite and affine Weyl groups, (see \cite{St} and \cite[Proposition 5.3.3]{MP},

\item{} simply laced extended affine Weyl
groups of rank $>1$ (see \cite[Theorem, III.1.14]{K}),

\item{} extended affine Weyl groups of index
zero, including extended affine Weyl groups of types $F_4$ and $G_2$ (see \cite[Theorem, 5.15]{A4}),

\item{} nullity $2$ extended affine Weyl groups of types $A_1$, $B_\ell$, $C_\ell$ (see \cite[Theorem 3.1]{A3}),

\item{} all but one, nullity $3$  extended affine Weyl groups of type $A_1$ (see \cite[Theorem 5.16 and Corollary 5.17]{AS2}).
\end{itemize}
It is shown in \cite{AS2}, that for each
 nullity $>2$, there is at least one extended affine Weyl group of type $A_1$
 which does not have the presentation by conjugation.
In \cite{H1} a necessary and sufficient condition is given for the existence of the presentation by conjugation for an extended affine Weyl group (minimality of the corresponding root system). Also in \cite{AS2},
using a notion of minimality on the set of generators of the Weyl group, a necessary and sufficient condition is given for the existence of the presentation by conjugation for $A_1$-type extended affine Weyl groups. One of the most interesting results related to this presentation is that the kernel of the epimorphism $\psi$ defined above, in the case of $A_1$, is a
direct sum of a finitely many copies of $\bbbz_2$,
where an upper bound is found for the number of copies. We show that this result remains valid for all reduced extended affine Weyl groups (Corollary \ref{fgag}).

In this work we consider non-simply laced extended affine root systems (Weyl groups). Since the problem of
existence of the presentation by conjugation for
extended affine Weyl groups of type $BC_\ell$ can
be reduced to those of type $B_\ell$ (see [Lemma 5.7]\cite{H1}), we only concentrate on types $B_\ell$, $C_\ell$, $F_4$ and $G_2$.

In Section \ref{preliminaries}, we record some basic facts regarding semilattices which will be of our use in the forthcoming sections. We have tried to keep this section as short as possible, we refer the interested reader to
\cite{AABGP}, \cite{A4} and \cite{AS2} for details on the topics involved.

In Section \ref{presentation}, we associate a presented group
$\hat{\w}$ to the extended affine root system $R$ (see Definition \ref{DEF}). First, by using only the defining relations of $\hat{\w}$, we find a suitable finite set of generators for $\hat\w$ and its center and also determine the relations among these generators (see Lemmas \ref{induces}, \ref{conj-cen-1}--\ref{well-define} and Proposition \ref{pox}(i)-(vii)). Next, we show that
kernel of $\psi$ is isomorphic to the direct sum of $n_0=\log_2 n$ copies of $\bbbz_2$,
where  $n$ is the number of integral collections for $R$ (see  Proposition \ref{ro}).
Finally, we state our main
theorem (Theorem \ref{Theorem-3}) which reveals several necessary and sufficient conditions for the existence of the presentation by conjugation for a given extended affine Weyl group. In particular, $\w$ has the presentation by conjugation if and only if $R$ has only one integral collection, namely the trivial collection. This provides a computational tool by which one can determine whether or not a given extended affine Weyl group has the presentation by conjugation.

As an application of the main theorem,
we determine a large class of extended affine Weyl groups which have (or don't have) the presentation by conjugation. In particular, up to
isomorphism, we give the precise answer for the existence of the presentation by conjugation for all extended affine Weyl groups of nullities $\leq 3$ (see Corollaries \ref{Cor-B2}--\ref{exist}). We concluded the paper with an appendix which provides a second proof for one of the earlier results in the paper.

The authors would like to thank A. Abdollahi for a fruitful discussion on Hopfian groups which led to the proof of Lemma \ref{AA}.

\vspace{5mm}
\section{\bf Preliminaries}\setcounter{equation}{0}\label{preliminaries}
\subsection{\bf
Essential supporting class of a
semilattice}
\setcounter{equation}{0}\label{semi}

For any positive integer $n$, we write $J_n=\{1,\ldots,n\}$ and
for $t\in J_n$, we set $J^t_n=\{t+1,\ldots,n\}$. Also if $r,s\in
J_n$ with $r<s$, we write $r<s\in J_n$, in this case also we
denote $\{r,s\}$ with $\{r<s\}$. We use a similar notation for
$r\leq s$. If $J\subseteq J_n$ and $x_r$'s are elements of a
group $G$, by $\prod_{r\in J}x_r$ we mean the product with the
usual order on $J$ as a subset of $J_n$. If $J$ is an empty set, we interpret the product as zero.

Let $S$ be a {\it semilattice} in a $\nu$-dimensional real vector
space $\v^0$, that is, $S$ is a discrete spanning subset of $\v^0$
satisfying $0\in S$ and $S=S\pm 2S.$ The $\bbbz$-span $\la S\ra$
of $S$ is a free abelian group of rank $\nu$, called a {\it
lattice} in $\v^0$. By \cite[Proposition II.1.11]{AABGP}, there
exists a set $\matB=\{\sg_1,\ldots,\sg_\nu\}$ satisfying
\begin{equation}\label{basis-1}
\matB\sub S,\quad\v^0=\sum_{r=1}^\nu\bbbr\sg_r\andd\la
S\ra=\sum_{r=1}^{\nu}\bbbz\sg_r. \end{equation}  We call such a
set $\matB$, a {\it basis}\index{basis} for $S$ and we fix it
throughout this work. For a set $J\sub J_\nu$ we put
\begin{equation}\label{sumo}
\tau_{_J}:=\sum_{r\in J}\sg_r\in \la S\ra
 \end{equation} (If $J=\emptyset$ we have by convention $\sum_{r\in
J}\sg_r=0$). With respect to $\matB$ we define
\begin{equation}\label{defsup}\supp_\matB(S)=\big\{J\sub
J_\nu \;:\; \tau_{_J}\in S\big\}
\end{equation}
\andd
\begin{equation}\label{defesup}
\esupp_\matB(S)=\big\{J\in \supp_\matB(S) \;:\; |J|\geq3\big\}.
\end{equation} Since $S\pm 2S\sub S$, it follows that
$$S=\biguplus_{J\in\supp_\matB(S)}(\tJ+2\la S\ra).
$$
 The collections $\supp_\matB(S)$ and $\esupp_\matB(S)$ are called
the {\it supporting class} and {\it essential supporting class} of
$S$ (with respect to the basis $\matB$), respectively. We drop
the subscript $\B$ when there is no confusion. The integer
$|\supp_\matB(S)|-1$ is called the {\it index} of $S$ and is
denoted by $\ind(S)$. Since $\matB\sub S$, we have
$\{\{\},\{1\},\ldots,\{\nu\}\}\sub\supp_\matB(S)$, and so
$\nu\leq\ind(S)\leq 2^\nu$.

\subsection{\bf Extended affine Weyl groups}\label{extended affine}

Let $\v$ be a finite dimensional real vector space and $\fm$ be a
symmetric positive semidefinite bilinear form on $\v$. An element
$\a\in\v$ is called {\it non-isotropic} (resp. {\it isotropic}) if
$(\a,\a)\not=0$ (resp. $(\a,\a)=0$). For a non-isotropic element $\a$,
we set $\ac=2\a/(\a,\a)$. The set of non-isotropic elements of a
subset $A$ will be denoted by $A^\times$.

Throughout this work we
assume $(\v,\fm,R)$ is a reduced non-simply laced extended affine
root system of rank $\ell$, nullity $\nu$ and twist
number $t$ (see \cite[Chapter II]{AABGP} for details). As it was explained in the Section \ref{introduction}, the results for type $BC_\ell$ can be reduced to those of type $B_\ell$. We denote the type of $R$ with $X$. So
$X=B_\ell (\ell\geq 2),\; C_\ell (\ell\geq
3),\;F_4$, or $G_2$. Let $\v^0$ be the radical of the form. By \cite[Lemma
II.4.15 and Proposition II.4.17]{AABGP}, we may find two subspaces
$\v^0_1$ and $\v^0_2$ of dimension $t$ and $\nu-t$, respectively,
and two semilattices $S_1$ and $S_2$ in $\v^0_1$ and $\v^0_2$,
respectively,  with $\vz=\v^0_1\oplus\v^0_2$ and
\begin{equation}\label{AABGP}
R=R(X, S_1, S_2):=(S+S)\cup(\rds+S_1\oplus\la
S_2\ra)\cup(\rdl+k\la S_1\ra\oplus S_2),
\end{equation}
where $S:=S_1\oplus\la S_2\ra$, and $\rd:=\{0\}\cup\rds\cup\rdl$
is an irreducible finite root system of type $X$ and rank $\ell$
with $\rds$ as the set of short roots and $\rdl$ as the set of
long roots of $\rd$. Also $k=3$ if $X=G_2$ and $k=2$, otherwise.
 Note that if $t=\nu$, then  $\vz_2=S_2=\{0\}$. By  \cite[Proposition II.4.9]{AABGP},
  if $X=F_4$ or $G_2$, then $S_1$ and $S_2$ are lattices in
$\vz_1$ and $\vz_2$, respectively. Also, if $X=B_\ell$ (resp.
$X=C_\ell$) with $\ell\geq3$, then $S_2$ (resp. $S_1$) is a
 lattice in $\vz_2$ (resp.  $\vz_1$). It is known from the classification
 of finite root systems that we may fix (we do) a basis
\begin{equation}\label{dotPi}
\dot\Pi=\{\a_1,\ldots,\a_\ell\}\hbox{ of }\rd\;\;\hbox{
with}\;\;\a_1\in\rds,\a_2\in\rdl\;\;\hbox{ and }\;\;
(\a_1,\a_2)\not=0.
\end{equation}
To emphasize on the roles of $\a_1$ and $\a_2$ in our work and to
distinguish them from other simple roots in $\dot{\Pi}$, we set
$$
\theta_1:=\a_1\andd\theta_2:=\a_2.$$
As in Subsection \ref{semi}, we may
fix a basis $\B=\B_1\cup\B_2$ of $S$ with
$\matB_1:=\{\sg_1,\ldots,\sg_t\}$,
$\matB_2:=\{\sg_{t+1},\ldots,\sg_\nu\}$ and
\begin{eqnarray}\label{sig} \matB_1\sub
S_1\quad\hbox{with}\quad\la S_1\ra=\sum_{r=1}^{t}\bbbz\sg_r\andd
\matB_2\sub S_2\quad\hbox{with}\quad\la
S_2\ra=\sum_{r=t+1}^{\nu}\bbbz\sg_r,
\end{eqnarray}
that is $\B_1$ (resp. $\B_2$) is a basis of $S_1$ (resp. $S_2$).

Here we need to recall some terms from \cite{AS3} which will be needed in the sequel. For any $i\in J_\ell$ and $r\in J_\nu$,  we set
\begin{equation}\label{t-ir}
 k_{i,r}:=\min\{n\in\mathbb N\mid \a_i+n\sg_r\in R\}.
  \end{equation}
   Then from (\ref{AABGP}) and    (\ref{sig})      we
 have
  \begin{equation}\label{fact-1}
\a_i+\bbbz k_{i,r}\sg_r\sub R.\quad\hbox{with}\quad
k_{i,r}=\left\{\begin{array}{ll}
   1, &\hbox{if}\;\a_i\in\rds,\vspace{2mm} \\
 k_r, &\hbox{if}\;\a_i\in\rdl,
   \end{array}\right.
\end{equation}
where
 \begin{equation}\label{kr}
k_r:=\left\{\begin{array}{ll}
k,&\hbox{if }r\in J_t,\vspace{2mm}\\
1,&\hbox{if }r\in J^t_\nu.
\end{array}\right.
  \end{equation}

Next, for all $ r\leq s\in J_\nu$ and  $i,j\in J_\ell$, we set
 \begin{eqnarray}\label{aij}
 a_{i,j}(r):= k_{j,r}k_{i,r}^{-1}(\a_i, \a_j^\vee),\qquad a_{i,j}(r,s):=
kk_r^{-1}k_{i,r}k_{j,s}(\a^\vee_i, \a_j^\vee),
   \end{eqnarray}
   \andd
\begin{eqnarray}\label{delta}
 \Delta(r,s)\hspace{-.5mm}:=\hspace{-.5mm}\left\{\begin{array}{ll}
 \hspace{-2mm}\d_1(r,s),      &\hbox{if\; $r<s\in J_t$},\vspace{2mm} \\
\hspace{-2mm}1,              & \hbox{ if $(r,s)\in J_t\times J^t_\nu$,}\vspace{2mm} \\
\hspace{-2mm}\d_2(r,s),       &\hbox{if $r<s\in J^t_\nu$,}
    \end{array}\right.
   \end{eqnarray}
   where\begin{eqnarray}\label{supofu}
 \d_j(r,s)\hspace{-.5mm}:=\hspace{-.5mm}\left\{\begin{array}{ll}
 \hspace{-2mm}1,& \hbox{if }\{r,s\}\in\supp(S_j),\vspace{2mm}\\
\hspace{-2mm} 2, &\hbox{if }\{r,s\}\not\in\supp(S_j).
  \end{array}\right.
   \end{eqnarray}
   By \cite[(2.18) and (2.19)]{AS3}, we have
\begin{eqnarray}\label{bvc}
 a_{i,j}(r)\in\bbbz\andd \Delta(r,s)^{-1}a_{i,j}(r,s)\in\bbbz.
   \end{eqnarray}

We now briefly recall the definition of an extended affine Weyl group. Let $\vd$ be the real span of $\rd$ and set $\vt:=\v\oplus(\v^0)^{\star}=\vd\oplus\v^0\oplus(\v^0)^\star$, where $(\v^0)^*$ is the
dual space of $\v^0$. Extend the form on $\v$ to $\vt$ naturally, by
dual pairing, namely
\begin{equation}\label{dual-pairing}
(\dot{\b}_1+\d_1+\lam_1,\dot{\b}_2+\d_2+\lam_2):=(\dot{\b}_1,\dot{\b}_2)+\lam_1(\d_2)+\lam_2(\d_1),
\end{equation} for $\dot{\b}_i\in\vd$, $\d_i\in\v^0$ and $\lam_i\in (\v^0)^*$. The {\it
(extended affine) Weyl group} $\w$ of $\r$ is by definition the
subgroup of $\GL(\vt)$ generated by reflections $w_\a$,
$\a\in\rtimes$, defined by $w_\a(u)=u-(u,\ac)\a$, $u\in\vt$. We
may identify the finite Weyl group $\dot{\w}$ of $\rd$ as a
subgroup of $\w$. We note that the following relations hold in
$\w$.
\begin{equation}\label{normal}
w_\a^2=1\andd ww_{\a}w^{-1}=w_{w(\a)}\qquad (\a\in R^\times,
w\in\w).
\end{equation}
Finally, we recall from \cite[Lemma 3.18(i)]{A4} or \cite[Proposition
2.1 (vii)-(viii)]{AS3} that
\begin{equation}\label{free-4}
\mbox{the center $Z(\w)$ of $\w$ is a free abelian group of rank $\nu(\nu-1)/2$.}
\end{equation}

We recall that a group $G$ is called {\it Hopfian}
if any epimorphism from $G$ onto $G$ is one-to-one. A group $G$ is said to satisfy Max-$n$ condition if
any ascending chain of its normal subgroups terminates. Finite groups and free abelian groups
 of finite rank satisfy Max-$n$ condition. Also if $N$ is normal in $G$ and both $N$ and $G/N$
 satisfy Max-$n$ condition then so does $G$. Finally, it is known  that any group satisfying Max-$n$ condition is
Hopfian (see \cite{R}, page 40).

\begin{lem}\label{AA}
$\w$ is a Hopfian group.
\end{lem}
\proof From Lemma 3.18 and \cite[Propositions 3.25, 3.27]{A4}, we know that $\w$ contains a normal subgroup $H$,
satisfying $\w\cong\dot{\w}\ltimes H$, where $Z(H)$ and
$H/Z(H)$ are free abelian groups of finite
rank. Since $Z(H)$, $H/Z(H)$ and $\w/H$ satisfy Max-$n$, $\w$ also does, and so is
 Hopfian.
\qed

\begin{rem}\label{remnew}
(i) The fact that $\w$ is Hopfian also follows from Theorem \ref{Theorem-3}.

(ii) If $R$ is the root system of an extended affine Lie algebra, then the Weyl group of $R$ which we defined here is isomorphic to
the Weyl group of the corresponding Lie algebra which, as usual, is defined as a
subgroup of the general linear group of the corresponding Cartan
subalgebra.
\end{rem}

 \section{\bf PRESENTATION BY
 CONJUGATION}\setcounter{equation}{0}\label{presentation}
 We keep all the notations and assumptions  as in the previous
section. In particular, $R$ is a reduced non-simply laced extended
affine root system of rank $\ell$,  nullity $\nu$ and twist number
$t$ of the form (\ref{AABGP}) and $\w$ is its extended affine Weyl
group. Throughout this work we denote the center of a group $G$ with $Z(G)$ and the commutator $x^{-1}y^{-1}xy$ of $x,y\in G$ with $[x,y]$. Recall that if $x,y\in Z(G)$, then
$xy=yx$, and so $(xy)^n=x^ny^n$ for all $n\in\bbbz$. Also if
$x,y,z\in G$ and $[x,y]\in Z(G)$, then $[x,yz]=[x,y][x,z]$ and $[x^n,z^m]=[x,z]^{nm}$ for all
 $n,m\in\bbbz$. These facts will be used in the sequel without any further reference.

\begin{DEF}\label{DEF} Let $\hat{\w}$ be the group defined by generators
$\hw_{\a}$, $\a\in\r^\times$ and relations:
$$\begin{array}{l}
  (\bi) \hspace{2mm}\hw_\a^2=1,\qquad \a\in\r^\times\vspace{2mm} \\
  (\bii)\hspace{2mm}\hw_\a\hw_\b\hw_\a=\hw_{w_{\a}(\b)},\qquad
   \a,\b\in\r^\times. \\
\end{array}$$
We say that the extended
affine Weyl group $\w$ of $R$ has {\it the presentation by
conjugation} if $\w\cong\hat{\w}$. It follows from Lemma \ref{AA} that, $\w$ has the presentation by conjugation if and only if
the epimorphism
\begin{equation}\label{main}
 \psi:\hat{\w}\longrightarrow\w,
 \end{equation}
induced by the assignment $\hw_{\a}\longmapsto w_{\a}$ is one-to-one. (Note that by (\ref{normal}), the defining relations
of $\hat{\w}$ are satisfied in $\w$.)
\end{DEF}


%

Since the finite Weyl group $\dot\w$ has the presentation by
conjugation (see [St]), the restriction of $\psi$ to
$\hat{\dot\w}:=\la \hw_{\a}\mid \a\in\dot{R}^\times\ra$ induces
the isomorphism
\begin{equation}\label{finite-case}
 \hat{\dot\w}\cong^{^{\hspace{-2mm}\psi}}\dot\w.
 \end{equation}
 One can easily deduce from relation (\bi) and (\bii) that
\begin{equation}\label{main1}
\hw\hw_{\a}\hw^{-1}=\hw_{\psi(\hw)(\a)}\qquad
(\hw\in\hat\w,\;\a\in\r^\times).
\end{equation}
\begin{lem}\label{induces}
(i) $\ker(\psi)\subseteq Z(\hat\w)$.

 (ii) $Z(\hat\w)=\psi^{-1}(Z(\w))$ and $\psi(Z(\hat\w))=Z(\w)$.

 (iii) $Z(\hat{\w})/\ker(\psi)\cong Z(\w)$.
\end{lem}
\proof (i) Let $\hw\in \ker(\psi)$. From (\ref{main1}) we have
$\hw\hw_\a\hw^{-1}=\hw_{\psi(\hw)(\a)}=\hw_\a$ for all $\a\in
R^\times$ and so $\hw\in Z(\hat\w)$. Thus (i) holds.

(ii) Let $\hw\in Z(\hat\w)$ and $\hz\in \psi^{-1}(Z(\w)))$. Then
From (\ref{normal}) and the fact that $\psi(\hz)(\a)=\a$ for all
$\a\in \v$ (see \cite[Proposition 2.1(vi), (2.13) and (2.7)]{AS3})
we have
$\psi(\hw)w_\a\psi(\hw)^{-1}=\psi(\hw\hw_\a\hw^{-1})=\psi(\hw_\a)=w_\a$
and $\hz \hw_\a \hz^{-1}=\hw_{\psi(\hz)(\a)}=\hw_{\a}$ for all
$\a\in R^\times$ and so $\hw\in \psi^{-1}(Z(\w)))$ and $\hz\in
Z(\hat\w)$. Thus (ii) holds. By (ii),
 the restriction of $\psi$ to $Z(\hat{\w})$ induces
the epimorphism $\psi:Z(\hat{\w})\longrightarrow Z(\w)$ and so
(iii) holds. \qed

Let $\supp(S_j)$,  $j=1,2$, be the supporting class
 of $S_j$ with respect to
$\B_j$ (see (\ref{defsup})) and $\dot\Pi$ be the basis of $\rd$ in
the form (\ref{dotPi}). We set
\begin{equation}\label{Pi}
\Pi:=\dot\Pi\cup\Pi_X, \qquad\hbox{where}
\end{equation}
  \begin{equation*}
\Pi_X:=\hspace{-1mm}\left\{\begin{array}{ll}
  \hspace{-2mm}\big\{\theta_1+\tJ\mid J\in
  \supp(S_1)\}\cup\big\{\theta_2+\tJ\mid J\in
\supp(S_2)\},&\hbox{if $X=B_2$,}\vspace{2mm}\\
\hspace{-2mm} \big\{\theta_1+\tJ\mid J\in\supp(S_1)\big\}\cup
   \big\{\theta_2+\sg_r\mid r\in J^t_\nu\big\}, &\hbox{if $X=B_{\ell}(\ell\geq3)$,}\vspace{2mm}\\
   \hspace{-2mm} \big\{\theta_1+\sg_r\mid r\in J_t\big\}\cup
   \{\theta_2+\tJ\mid J\in\supp(S_2)\big\}, &\hbox{if $X=C_{\ell}(\ell\geq3)$,}\vspace{2mm}\\
   \hspace{-2mm}\big\{\theta_1+\sg_r,\theta_2+\sg_s\mid r\in J_t,\; s\in J^t_\nu\big\},&
  \hbox{if $X=F_4$ or $G_2$.}
    \end{array}\right.\end{equation*}

An argument similar to \cite[Section 4]{A4} gives the following
lemma (see \cite[Theorem 2.3.3]{Sh} for a detailed proof).
    \begin{lem}\label{basis}
$\w=\la w_\a:\a\in\Pi\ra$ and $\w\Pi=\rtimes$.
 \end{lem}
From the way $\hat{\w}$ is defined we have $\hat{\w}=F/N$, where
$F$ is the free group on the set $\{r_\a\mid\a\in \r^\times\}$ and
 $N$ is the normal closure of the set $K:=\{r_\a^2,\;
 r_\a r_\gamma r_\a r_{w_\a(\gamma)}\mid\a,\gamma\in\r^\times\}$
in $F$. Then $\hw_{\a}=r_\a N$, $\a\in\rtimes$. Now fix
$\beta\in\rtimes$ and let $\w\b$ be the orbit of $\b$ under $\w$.
Let $\Phi_\b:F\longrightarrow\bbbz_2$ be the epimorphism induced
by the assignment
 $$\Phi_\b(r_{\a})=\left\{\begin{array}{ll}
   1, & \hbox{if}\;\; \a\in\w\b,\vspace{2mm} \\
   0, & \hbox{if}\;\; \a\in\r^\times\setminus\w\b.
 \end{array}\right.$$
Since $K\sub\ker(\Phi_\b)$, the epimorphism $\Phi_\b$ induces a
unique epimorphism $\Phi_\b:\hat{\w}\longrightarrow\bbbz_2$ so
that $\bar{\Phi}_\b(\hw_\a)=\Phi_\b(r_\a)$, $\a\in\r^\times$.
  Then
  \begin{equation}\label{subset}
   \bar{\Phi}_\b(\la \hw_\a\mid
   \a\in\r^\times\setminus\w\b\ra)=\{0\}\andd
   \bar{\Phi}_\b(\hw_{\b_1}\cdots\hw_{\b_n})=\bar{n},\quad(\b_i\in\w\b),
\end{equation}
where $\bar{n}$ is the image of $n$ in $\bbbz_2$ under the natural
map.
\begin{lem}\label{minimal} $\{\hw_{\a}\mid\a\in\Pi\}$
 is a minimal set of generators for $\hat{\w}$.
\end{lem}

\proof If  $\a\in\rtimes$, then by Lemma \ref{basis} we have
$\a=w_{\b_1}\cdots w_{\b_m}(\b)$ for some $\b_i\in\Pi$ and so by
(\ref{main1}),
$\hw_\a=\hw\hw_{\b}\hw^{-1}\in\la\hw_\a\mid\a\in\Pi\ra$ where
$\hw=\hw_{\b_1}\cdots\hw_{\b_m}$. Thus the set in the statement
generates $\hat\w$.  To show that
 it is minimal, fix $\b\in\Pi$ and set
$\Pi_\b:=\Pi\setminus\{\b\}$. We must show that the elements
$\hw_\a$, $\a\in\Pi_\b$ do not generate $\hat{\w}$. We show this
by a type dependent argument. First we note from \cite[Lemma
2.16]{A4} that
$$\w\b\subseteq\left\{\begin{array}{ll}
\b-\theta_1+\rds+\la L\ra &\hbox{if }X=B_\ell(\ell\geq 2),\;\b\in\rs,\\
\b-\theta_2+\rdl+2\la S\ra&\hbox{if }X=B_2\hbox{ or
}X=C_\ell(\ell\geq 3),\;\b\in\rl.
\end{array}\right.
$$
From this (and the way $\Pi$ is defined) one sees that
$\Pi_\b\sub\rtimes\setminus\w\b$.
 So form (\ref{subset}) we
have $\bar{\Phi}_\b(\hw)=0$ for any
  $\hw\in\la\hw_\a\mid\a\in\Pi_\b\ra$ and
 $\bar{\Phi}_\b(\hw_\b)=1$. Thus
 $\hw_\b\not\in\la\hw_\a\mid\a\in\Pi_\b\ra$.
Next let, $\b\in\rl$ for $X=B_\ell(\ell\geq 3)$ or $\b\in\rs$ for
$X=C_\ell(\ell\geq 3)$. Then form the way $\Pi$ is defined, it is
clear that $\Pi_\b$ spans a $(\nu+\ell-1)$-dimensional subspace of
$\v$. Now suppose to the contrary that $\hw_{\b}=\hw_{\b_1}\cdots
\hw_{\b_m}$ for some $\b_i\in\Pi_\b$. Then applying the
epimorphism $\psi$ we get $w_{\b}=w_{\b_1}\cdots w_{\b_m}$. By
acting both sides on $\b$ we see that
$\b\in\span_{\bbbr}\{\b_1,\ldots,\b_m\}\subseteq\span\Pi_\b$. That
is $\span\Pi_\b=\span\Pi=\nu+\ell$, a contradiction. The argument
for types $F_4$ and $G_2$ is exactly the same as what we did in
 this paragraph. This completes the proof.\qed

For $\a\in\rtimes$, $\sg\in\v^0$ with $\a+\sg\in R$, we set
\begin{equation}\label{deft}
\hat{t}^\sg_\a=\hw_{\a+\sg}\hw_{\a}.
 \end{equation}
 The elements
of this form play a crucial role in our description of $\hat{\w}$.

\begin{con}\label{convention} We reserve the symbols
$\a,\b,\gamma$ for elements of $\rtimes$ and symbols
$\delta,\sg,\eta$ for elements of $\v^0$. By a symbol such as
$\hat{w}_{\a+\sg}$ or $\hat{t}^\sg_\a$, we always mean that $\a\in\rtimes$,
$\sg\in\v^0$ and $\a+\sg\in R$. In places that this convention
might cause confusion we clarify the situation. To see how we use
this convention note for example that from Lemma \ref{conj-cen-1}(iii) we understand that $\d_r\in\v^0$ and
$\a,\a+\d_r,\a+\sum_{r=1}^{n}\d_r\in R$, for all $r$.
\end{con}

\begin{lem}\label{conj-cen-1} Any element in $\hat{\w}$ which has
 one of the following forms is central:
$$
\begin{array}{ll}
(i)&[\th_\a^\sg,\th_\b^\d],\vspace{2mm}\\
 (ii)&\th_{\a+\sg}^\d\th_\a^{-\d},\vspace{2mm}\\
(iii)&\th_\a^{-\Sigma_{r=1}^n\d_r}\prod _{r=1}^n\th_\a^{\d_r}.
\end{array}
$$
\end{lem}

\proof  Since the image of each element of the given forms is
central in $Z(\w)$ (see \cite[Lemmas 1.1(iii)-(vii) and 1.2]{AS1}
for details), the result follows immediately from Lemma
\ref{induces}(ii).\qed

Here we record the following useful lemma from \cite[Lemma 3.18]{A3}. For the convenience of the
reader we present its proof here with a simpler argument than the one in \cite{A3}. We recall
from \cite[Chapter II]{AABGP} that if $\a,\a+\sg\in\rtimes$ for
some $\sg\in\v^0$, then $\a+n\sg\in\rtimes$ for all $n\in\bbbz$.

\begin{lem}\label{gen-3} $(\th_\a^\sg)^n=
\th_\a^{n\sg}$ and $\th_{\a+n\sg}^\sg=\th_\a^\sg$ for all
$n\in\bbbz$. Moreover, $\th^{-\sg}_\a=\th^\sg_{-\a}$.
\end{lem}

\proof Clearly the first claim holds for $n=0,1$, moreover, for
$n\in\bbbz$ we have $w_{\a+\sg}w_{\a}
\big(\a+(n-2)\sg\big)=\a+n\sg$. Then using induction on $n\geq 0$
and applying relations (\bi), (\bii), we obtain
\begin{eqnarray*}
\th_\a^{n\sg}=\hw_{\a+n\sg}\hw_{\a}&=
&[(\hw_{\a+\sg}\hw_{\a})\hw_{\a+(n-2)\sg}(\hw_{\a+\sg}\hw_{\a})^{-1}]
\hw_{\a}\\
&=&(\hw_{\a+\sg}\hw_{\a})(\hw_{\a+(n-2)\sg}\hw_{\a})(\hw_{\a+\sg}
\hw_{\a})\\
&=&(\hw_{\a+\sg}\hw_{\a})(\hw_{\a+\sg}\hw_{\a})^{n-2}(\hw_{\a+\sg}\hw_{\a})\\
&=&(\hw_{\a+\sg}\hw_{\a})^n=(\th_\a^\sg)^n.
\end{eqnarray*}
If $n<0$, then $(\hw_{\a+\sg}\hw_{\a})^{n}=
(\hw_{\a-\sg}\hw_{\a})^{-n}=w_{\a+n\sg}\hw_{\a}$. The second claim
holds since
$$\th^\sg_{\a+n\sg}=\hw_{\a+n\sg+\sg}(\hw_\a\hw_\a)\hw_{\a+n\sg}=
\th^{(n+1)\sg}_\a\th^{-n\sg}_\a=\th_\a^\sg.
$$
Finally, since by (\ref{main1}),
$\hw_\a\hw_{\a+\sg}\hw_\a=\hw_{-\a+\sg}$, the last assertion
holds.
 \qed

\begin{lem}\label{ip} (i) $\th^{-\d}_{\a+\sg}\th^{\d}_\a=\th^\sg_{\a+\d}\th^{-\sg}_\a$.

(ii)
$[\th_\a^\sg,\;\th_\b^\d]=\th^{-\d}_{\b-(\b,\ac)\sg}\th^\d_\b=
\th^{-(\b,\ac)\sg}_{\b+\d}\th^{(\b,\ac)\sg}_\b$. In particular, if
$(\a,\b)=0$ or $\sg=\d$, then $[\th_\a^\sg,\;\th_\b^\d]=1$.
\end{lem}

\proof (i) By Lemmas \ref{conj-cen-1}(ii) and \ref{gen-3}, we have
$$\th^{-\d}_{\a+\sg}\th^{\d}_\a=\hw_{\a+\sg}\hw_{\a+\sg+\d}\hw_{\a+\d}\hw_\a=
\hw_{\a+\sg+\d}\hw_{\a+\d}\hw_\a\hw_{\a+\sg}=\th^\sg_{\a+\d}\th^{-\sg}_\a.
$$
(We have used the fact that if $xy$ is central for two elements
$x$, $y$ of a group $G$, then $xy=yx$.)

(ii) Using (\ref{main1}) and Lemma \ref{gen-3}, we have
$[\th_\a^\sg,\;\th_\b^\d]=\th^{-\sg}_\a\th^{-\d}_\b\th^\sg_\a\th^\d_\b
=\th^{-\d}_{\psi(\th^{-\sg}_\a)(\b)}\th^\d_\b.$ But
$\psi(\th^{-\sg}_\a(\b))=w_{\a}w_{\a+\sg}(\b)=\b-(\b,\ac)\sg$. The
second equality follows immediately from part (i) and (\ref{deft}). Finally if
$(\a,\b)=0$ or $\sg=\d$, then using Lemma \ref{gen-3}, we are
done.$\quad$ \qed

\begin{lem}\label{well-define}
(a) If $\a,\b$ belong to the same orbit of $\rtimes$, under the
action of $\w$, then
$$
\begin{array}{ll}
(i)&[\th_{\a}^\sg,\th_{\a}^\d]=
 [\th_{\b}^\sg,\th_{\b}^\d],\vspace{2mm}\\
 (ii)&\th_{\a+\sg}^{\d}\th_{\a}^{-\d}=\th_{\b+\sg}^{\d}
\th_{\b}^{-\d},\vspace{2mm}\\
(iii)&\th_{\a}^{-\Sigma_{r=1}^n\d_r}\prod _{r=1}^n\th_{\a}^{\d_r}=
\th_{\b}^{-\Sigma_{r=1}^n\d_r}\prod _{r=1}^n\th_{\b}^{\d_r}.
\end{array}
$$

(b) If $\a,\a'$ belong to the same orbit of $\rtimes$
 under the action of $\w$ and $\b,\b'$ are elements of $\rtimes$ such that $(\a,\bc)=(\a',{\b'}^\vee)$ then
 $$[\th_{\b}^{\d},\th_{\a}^{\sg}]=
 [\th_{\b'}^{\d},\th_{\a'}^{\sg}].$$
\end{lem}

\proof (a) Let $w\in\w$ be such that $\b=w(\a)$ and fix a primage
$\hw\in\hat{\w}$ of $w$ under $\psi$. By Lemma \ref{conj-cen-1},
the left hand sides (LHS) of the equalities in the statement are all
central, and so $\hw\hbox{LHS}\hw^{-1}=\hbox{LHS}$. But by
(\ref{main1}), we have $\hw\hw_{\a+\eta}\hw^{-1}=\hw_{\b+\eta}$
for any $\eta\in\v^0$. The result now follows immediately.

(b) Let $w\in\w$ be such that $\a'=w(\a)$ and fix a primage
$\hw\in\hat{\w}$ of $w$ under $\psi$.  Then using Lemma
\ref{ip}(ii) and (\ref{main1}) we have
 \begin{eqnarray*}
[\th_{\b}^{\d},\th_{\a}^{\sg}]=\hw[\th_{\b}^{\d},\th_{\a}^{\sg}]\hw^{-1}=
\th^{-\sg}_{w(\a)-(\a,\b^\vee)\d}\th^\sg_{w(\a)}=
\th^{-\sg}_{\a'-(\a',\b'^\vee)\d}\th^\sg_{\a'}=[\th_{\b'}^{\d},\th_{\a'}^{\sg}].
\end{eqnarray*}
This completed the proof.\qed

For further study of the center of $\hat{\w}$ we need to introduce some new terms. Our motivation for defining the term $\hat{z}_{_J}$ below has been \cite[Lemma 2.3]{AS3}. For $J\subseteq J_\nu$, $\a,\a'\in\rds$ and
$\b,\b'\in\rdl$ with $(\a',\b')<0$, we set
\begin{equation}\label{ZJ}
\hspace{-.17cm}\hz_{_J}:=\left\{\begin{array}{ll}
 \th_{\a}^{-\tJ}\prod_{r\in
    J}\th_{\a}^{\sg_r},& \hbox{if $J\in\mbox{supp}(S_1),$}\vspace{3mm}\\
  \big[\th_{\a}^{\sg_r},\th_{\a}^{\sg_s}\big], &  \hbox{if
   $J=\{r,s\}\not\in\mbox{supp}(S_1)$,
    $r<s\in J_t,$}\vspace{3mm}\\
    \th_{\b}^{-\tJ}\prod_{r\in
    J}\th_{\b}^{\sg_r},& \hbox{if $J\in\mbox{supp}(S_2),$}\vspace{3mm}\\
     \big[\th_{\b}^{\sg_r},\th_{\b}^{\sg_s}\big],&
     \hbox{if $J=\{r,s\}\not\in\mbox{supp}(S_2)$, $r<s\in J^t_\nu$,}\vspace{3mm}\\
 \big[\th_{\b'}^{\sg_s},\th_{\a'}^{\sg_r}\big],  & \hbox{if $J=\{r,s\}$,
 $r\in J_t,\; s\in J^t_\nu$ ,}\vspace{3mm}\\
    1, & \hbox{otherwise}.
\end{array}\right.\hspace{-.2cm}
\end{equation}
 We note that by
Lemma \ref{well-define}, definition of $\hz_{_J}$ does not depend
on the particular choice of $(\a,\b)\in\rds\times\rdl$ and
$(\a',\b')\in\rds\times\rdl$ with $(\a',\b')<0$ \big(if
$\a',\a''\in\rds,\b',\b''\in\rdl$ with $(\a',\b')<0$ and
$(\a'',\b'')<0$, then $(\a',\b'^\vee)=(\a'',\b''^\vee)$\big). Note
also that $\hz_{\{r\}}=1$ for all $r$. By (\ref{fact-1}),
(\ref{deft}) and Lemmas \ref{conj-cen-1} and \ref{induces}(ii) we
 have (for $J\subseteq J_\nu$ and $(i, r)\in J_\ell\times
J_\nu$)
 \begin{equation}\label{tir}
  \th_{i,r}:=\th^{k_{i,r}\sg_r}_{\a_i}=\hw_{\a_i+k_{i,r}\sg_r}\hw_{\a_i} \andd \hz_{_J}\in
  Z(\hat\w).
\end{equation}

The following proposition plays an essential role in the
sequel. Our argument for part (iii) of this proposition relies on several
known results in the literature. In particular, a known fact that
the presentation by conjugation holds for all extended affine
Weyl groups of nullities $\leq 2$ is used in the argument. For
convenience of the reader we give another proof for
Proposition \ref{pox}(iii) in the appendix (see Lemma
\ref{impor}) which assumes no previous knowledge
about the presentation by conjugation for low nullities.

Let us set
\begin{equation}\label{JX} \mathcal J=\mathcal J(X, S_1,
S_2):=\left\{\begin{array}{ll}
  \esupp(S_1)\cup\esupp(S_2) & \hbox{if $X=B_2$,}\vspace{2mm} \\
   \esupp(S_1) & \hbox{if $X=B_\ell(\ell>2)$,}\vspace{2mm} \\
   \esupp(S_2) & \hbox{if $X=C_\ell(\ell>2)$,}\vspace{2mm} \\
   \emptyset & \hbox{if $X=F_4$ or $G_2$.}
 \end{array}\right.
  \end{equation}
($\esupp(S_j)$, $j=1,2$  is given by  (\ref{defesup}).)

\begin{pro}\label{pox}
{(i)}
$\hat{\w}=\la\hw_{\a_i},\;\hat{t}_{i,r},\;\hz_{_{\{r,s\}}},\;\hz_{_J}\mid
i\in J_\ell,\; r\leq s\in J_\nu,\;J\in\J\ra$,

{(ii)} If $i,j\in J_\ell$,  $r\in J_\nu$ and $a_{i,j}(r)$ is
given by {(\ref{aij})}, then
$$\hw_{\a_i}\th_{j,r}\hw_{\a_{i}}=\th_{j,r}\th_{i,r}^{
-a_{i,j}(r)},$$

{(iii)} If $i,j\in J_\ell$,  $r\leq s\in J_\nu$, then
$$[\hat{t}_{i,r},\hat{t}_{j,s}]=
\hz_{_{\{r,s\}}}^{\Delta(r,s)^{-1}a_{i,j}(r,s)},$$  where
$a_{i,j}(r, s)$ and $\Delta(r,s)$ are given by {(\ref{aij})}
and {(\ref{delta})},

{(iv)}  If $J\in\supp(S_1)\cup\supp(S_2)$, then
 $$\hz_{_{J}}^2=\prod_{\{r,s\in
J~:~r<s\}}\hz_{_{\{r,s\}}}^{2\Delta(r,s)^{-1}},$$

{(v)} $Z(\hat\w)=\la \hz_{_{\{r,s\}}},\;\hz_{_J}\mid r<s\in
J_\nu,\;J\in\J\ra$, where $\J$ is given by { (\ref{JX})}.
\end{pro}
\proof (i) Let $T$ be  a  group given  in the right hand side of
the statement. By (\ref{tir}) and Lemma \ref{minimal} we only need
to show that for any $\b\in\Pi$, $\hw_\b\in T$. Clearly this holds
if $\b=\a_i$ for some $i\in J_\ell$. If $\b=\theta_j+\sg_r$ for
some $(j,r)\in (\{1\}\times J_t)\cup (\{2\}\
  \times J_\nu^t)$, then from the fact that
$\hw_\b=\th_{j,r}\hw_{\theta_j}$, it is clear that  $\hw_\b\in T$.
Next, let $\b=\theta_j+\tJ$, for some $j\in\{1,2\}$ and
$J\in\supp(S_j)$, then from the way $\hz_{_J}$ is defined we have
$\hw_{\theta_j+\tJ}=\hw_{\theta_j}\hz_{_J}(\prod_{r\in
J}\th_{j,r})^{-1}\in T$.

(ii) We have
\begin{eqnarray*}
\hw_{\a_i}\th_{j,r}\hw_{\a_{i}}&=&\th_{j,r}\th_{j,r}^{-1}\hw_{\a_i}\th_{j,r}\hw_{\a_{i}}\\
\hbox{(using (\ref{main1}))}&=&\th_{j,r}\hw_{\psi(\th_{j,r}^{-1})(\a_i)}\hw_{\a_{i}}\\
\hbox{(using (\ref{main}), (\ref{tir}) and
\bi)}&=&\th_{j,r}\hw_{w_{\a_i}w_{\a_i+k_{i,r}\sg_r}(\a_i)}\hw_{\a_{i}}\\
\hbox{(using (\ref{aij}))}&=&\th_{j,r}\hw_{\a_i-k_{j,r}(\a_i,\a^\vee_j)\sg_r}\hw_{\a_{i}}=\th_{j,r}\hw_{\a_i-a_{i,j}(r)k_{i,r}\sg_r}\hw_{\a_{i}}\\
\hbox{(using (\ref{bvc}) and  Lemma
\ref{gen-3})}&=&\th_{j,r}(\hw_{\a_i+k_{i,r}\sg_r}\hw_{\a_{i}})^{-a_{i,j}(r)}=\th_{j,r}\th_{i,r}^{-a_{i,j}(r)}.
\end{eqnarray*}

(iii) Set $\hz:=[\hat{t}_{i,r},\hat{t}_{j,s}]
\hz_{_{\{r,s\}}}^{-\Delta(r,s)^{-1}a_{i,j}(r,s)}$.  From
\cite[Lemma 2.5(i)]{AS3}, it follows that $\psi(\hz)=1$. If
$(\a_i,\a_j)=0$, then by Lemma \ref{ip}(ii), $\hz=1$. So we may
assume $(\a_i,\a_j)\not=0$. Set
$$T=R^\times\cap(\bbbr\a_i+\bbbr\a_j+\bbbr\sg_r+\bbbr\sg_s).
$$
Clearly $T$ satisfies the axioms of a Saito extended affine root
system (\cite[Definition 1]{A2}) and so by \cite[Theorem 18]{A2},
$T$ is the set of non-isotropic roots of an extended affine root
system $R_T$ of nullity $\leq 2$ (see also \cite[Proposition
5.9]{H1}). Let $\hat{\w}_T$ be the group defined similar to
$\hat{\w}$, corresponding to the extended affine Weyl group of
$R_T$. From (\ref{tir}) and (\ref{ZJ}) and the way $R_T$ is
defined, it follows that $\hz\in \hat{\w}_T$. By \cite[Theorem 3.1]{A3}, the
Weyl group $\w_T$ has the presentation by conjugation, and so the
restricted map $\psi: \hat{\w}_T\longrightarrow\w_T$ is an
isomorphism.  Thus $\hz=1$.

(iv) Let $J\in\supp(S_j)$, $j\in\{1,2\}$ and set $\hw=\prod_{r\in
J} \th_{j,r}$. Then $k_{j,r}=k_{j,s}=1$ for all $r,s\in J$. Since
$\psi(\hat{t}_{j,r})(\theta_j)=\theta_j+2\sg_r$ for all $r\in J$,
we have $\psi(\hw)(\theta_j)=\theta_j+2\tau_{_J}$.  So from
(\ref{ZJ}) and the fact that $a_{j,j}(r,s)=a_{j,j}(r)=2$ (see
(\ref{aij})) we have
\begin{eqnarray*}
\hz_{_J}^2\hspace{-3mm}&=&\hspace{-3mm}(\th_{\theta_j}^{-\tJ}\hw)^2\\
    \hbox{(using  Lemmas
\ref{conj-cen-1}(iii) and \ref{gen-3}
)}\hspace{-3mm}&=&\hspace{-3mm}(\th_{\theta_j}^{-\tJ})^2\hw^2=(\th_{\theta_j}^{-2\tJ})\hw^2\\
\hbox{(using (\ref{deft}) and (\ref{tir}))
}\hspace{-3mm}&=&\hspace{-3mm}\hw_{\theta_j}\hw_{\theta_j+2\tau_{_J}}
\hw^2=\hw_{\theta_j}\hw_{\psi(\hw)(\theta_j)}
\hw^2\\
\hbox{(using (\ref{main1}))}\hspace{-3mm}&=&\hspace{-3mm}
\hw_{\theta_j}\hw\hw_{\theta_j}\hw^{-1}
\hw^2=\hw_{\theta_j}\hw\hw_{\theta_j}\hw\\
\hbox{(using $\bi$  and
(i))}\hspace{-3mm}&=&\hspace{-3mm}\big(\prod_{r\in J}
\hw_{\theta_j}\th_{j,r}\hw_{\theta_j}\big)(\prod_{r\in J}
\th_{j,r})=(\prod_{r\in J}
\th_{j,r}^{-1})(\prod_{r\in J} \th_{j,r})\\
\hbox{(using (iii))}\hspace{-3mm}&=&\hspace{-3mm} \prod_{r<s\in
J}[\th_{j,r},\th_{j,s}]= \prod_{\{r,s\in
J~:~r<s\}}\hz_{_{\{r,s\}}}^{2\Delta(r,s)^{-1}}.
\end{eqnarray*}

(v) By (\ref{tir}) it is enough to prove $Z(\hat\w)$ is a subset
 of the right hand side of the equality. Let $\hat w\in Z(\hat\w)$.
 By parts (i)-(iii), we see that
  $\hat w$  can be written as an expression in
the form
 $ \hw=\hat{\dot w}
\prod_{r=1}^{\nu}\prod_{i=1}^{\ell}\th_{i,r}^{n_{i,r}}\hz,$
 where $\hat{\dot w}\in\hat{\dot{\w}}$,
   $ n_{i,r}\in\bbbz$ and $\hz\in \la \hz_{_{\{r,s\}}},\;\hz_{_J}\mid r<s\in
J_\nu,\;J\in\J\ra$.
    Then  we have
\begin{equation*}
  \psi(\hw)=\psi(\hat{\dot w})\prod_{r=1}^{\nu}
  \prod_{i=1}^{\ell}t_{i,r}^{n_{i,r}}\psi(\hz)\in Z(\w)\quad\mbox{where}\quad
   t_{i,r}:=\psi(\th_{i,r})=w_{\a_i+k_{i,r}\sg_r}w_{\a_i}.
  \end{equation*}
Since $\psi(\hat{\dot w})\in \dot{\w}$ and $\psi(\hz)\in Z(\w)$,
we have from \cite[(2.9) and Proposition 2.1(v)]{AS3} and
(\ref{finite-case}) that $\hat{\dot{w}}=1$ and $n_{i,r}=0$ for all
$i,r$. Thus $\hat{w}=\hat{z}\in\la \hz_{_{\{r,s\}}},\;\hz_{_J}\mid
r\leq s\in J_\nu,\;J\in\J\ra$. \qed

To obtain further information about $\ker(\psi)$ and consequently
about the existence of the presentation by conjugation for  extended affine Weyl groups under consideration, we need to introduce a new term,
called an {\it integral collection}.

For all $1\leq r<s\leq\nu$ and  $J\subseteq \{1,\ldots,\nu\}$, we
set
\begin{eqnarray}\label{chi}
 \chi_{_J}(r,s)=\left\{\begin{array}{ll}
 1,& \hbox{ if }\{r,s\}\subsetneq J,\vspace{2mm}\\
 0, &\hbox{otherwise,}
  \end{array}\right.
   \end{eqnarray}

  We   call
  $\epb=\{\ep_{_J}\}_{J\in\mathcal J}$, $\ep_{_J}\in\{0,1\}$,  an {\it integral collection}
  for $(S_1, S_2)$, if
  \begin{eqnarray}\label{integ-R}
 \epb_{r,s}:=\Delta(r,s)^{-1}\sum_{J\in\J}
 \chi_{_J}(r,s)\epsilon_{_J}\in\bbbz, \quad \hbox{ for all $r<s\in J_\nu$,}
       \end{eqnarray}
 where $\J$ is given by (\ref{JX}) (if $J$ is an empty set, we make the convention
that the sum is zero and interpret $\epb$ as the zero collection
). If $\ep_{_J}=0$ for all $J\in\J$, we call $\epb=\{0\}_{J\in\J}$
the {\it trivial} collection.  It is clear that the trivial
collection is an integral collection. Clearly there are at most
$2^{|\J|}$ integral collections for $(S_1, S_2)$. Any integral
collection different from the trivial collection is called {\it
non-trivial}.

For any integral collection $\epb=\{\ep_{_J}\}_{J\in\mathcal J}$
of $(S_1, S_2)$, we set
\begin{equation}\label{nbm2}
\hu(\epb):=\prod_{1\leq
r<s\leq\nu}\hz_{_{\{r,s\}}}^{\epb_{r,s}}\prod_{J\in\mathcal J }
\hz_{_{J}}^{\epsilon_{_J}}.
\end{equation}
 \begin{pro}\label{ro}
  {(i)} \begin{eqnarray*}
  \ker(\psi)&=&\big\{\hu\in Z(\hat\w) \mid \hu^2=1\big\}\\
 &=&\big\{\hu\in Z(\hat\w) \mid |\hu|<\infty\big\}\\
 &=&\big\{\hu(\epb) \mid \epb\mbox{ is an integral collection for }(S_1,
  S_2)\big\}.
 \end{eqnarray*}

{(ii)} The assignment $\epb\longmapsto\hu(\epb)$ is a one to
one correspondence from the set of integral collections for $(S_1,
S_2)$ onto $\ker(\psi)$.
    \end{pro}
   \proof (i) From Lemma \ref{induces}(iii) and (\ref{free-4}) we
have
\begin{eqnarray}\label{inq}
\big\{\hu\in Z(\hat\w) \mid \hu^2=1\big\}\subseteq\big\{\hu\in
Z(\hat\w) \mid |\hu|<\infty\big\}\subseteq\ker(\psi).
\end{eqnarray}
  Next, let $\hu\in\ker(\psi)$.  By Lemma \ref{induces}(i) and Proposition
\ref{pox}(iv)-(v), $\hu$  can be written in the form
\begin{equation}\label{nbm} \hu=\prod_{1\leq
r<s\leq\nu}\hz_{_{\{r,s\}}}^{m_{r,s}}\prod_{J\in\mathcal J }
\hz_{_{J}}^{\epsilon_{_J}},\qquad\mbox{($m_{r,s}\in\bbbz$,\;\;
$\epsilon_{_J}\in\{0,1\}$)}.
\end{equation}
 Then
 using  (\ref{nbm}),  (\ref{chi}) and Proposition \ref{pox}(v)   we have
    \begin{eqnarray*}
\hu^{2}&=&\prod_{1\leq
r<s\leq\nu}\hz_{_{\{r,s\}}}^{2m_{r,s}}\prod_{J\in\mathcal J }
\hz_{_{J}}^{2\epsilon_{_J}}\\
\mbox{(using Proposition \ref{pox}(iii))}&=&\prod_{1\leq
r<s\leq\nu}\hz_{_{\{r,s\}}}^{2m_{r,s}}\prod_{J\in\mathcal J}
(\prod_{\{r,s\in J~:~r<s\}}\hz_{_{\{r,s\}}}^{2\Delta(r,s)^{-1}})^{\epsilon_{_J}}\\
&=&\prod_{1\leq
r<s\leq\nu}\hz_{_{\{r,s\}}}^{2m_{r,s}}\prod_{J\in\mathcal J }
\prod_{1\leq
r<s\leq\nu}\hz_{_{\{r,s\}}}^{2\Delta(r,s)^{-1}\chi_{_J}(r,s)\epsilon_{_J}}\\
&=&\prod_{1\leq r<s\leq\nu}\hz_{_{\{r,s\}}}^{2m_{r,s}}
\prod_{1\leq
r<s\leq\nu}\hz_{_{\{r,s\}}}^{2\Delta(r,s)^{-1}\sum_{J\in\mathcal J
}\chi_{_J}(r,s)\epsilon_{_J}}
\end{eqnarray*}
\begin{eqnarray}\label{uhg}
\hspace{3cm}&=& \prod_{1\leq
r<s\leq\nu}\hz_{_{\{r,s\}}}^{2m_{r,s}+2\Delta(r,s)^{-1}\sum_{J\in\mathcal
J }\chi_{_J}(r,s)\epsilon_{_J}}
\end{eqnarray} and so
\begin{eqnarray}\label{uhg-1} 1=\psi(\hu^{2})= \prod_{1\leq
r<s\leq\nu}
z_{_{\{r,s\}}}^{2m_{r,s}+2\Delta(r,s)^{-1}\sum_{J\in\mathcal J
}\chi_{_J}(r,s)\epsilon_{_J}},
\end{eqnarray}
where $z_{_{\{r,s\}}}:=\psi(\hz_{_{\{r,s\}}})$. Set
$\epb:=\{\ep_{_J}\}_{J\in\mathcal J}$ and $\epb_{r,s}:=m_{r,s}$.
Since $\la z_{_{\{r,s\}}}\mid r<s\in J_\nu\ra$ is a free abelian
group on generators $z_{_{\{r,s\}}}$'s (see \cite[Lemmas 2.2(i)
and 2.3]{AS3}), the integers $\epb_{r,s}$'s satisfy
(\ref{integ-R}). Thus $\hu^2=1$ and $\hu=\hu(\epb)$ where $\epb$
is an integral collection. Thus
$$\ker(\psi)\subseteq\big\{\hu(\epb) \mid \mbox {$\epb$ is
 an integral collection for $(S_1, S_2)$}\big\}\subseteq \big\{\hu\in Z(\hat\w) \mid
\hu^2=1\big\}.$$ This together with (\ref{inq}) completes the
proof of (i).

 (ii) By (i), it is enough to show that the assignment $\epb\longmapsto\hu(\epb)$ is
 injective. Now let $\hu(\epb)=\hu(\epb')$, where  $\epb=\{\ep_{_J}\}_{J\in\J}$
  and $\epb'=\{\ep_{_J}'\}_{J\in\J}$ are two integral collections for $(S_1, S_2)$. We claim that
$\epb=\epb'$.  Suppose to the contrary that
$\ep_{_{J_0}}\neq\ep'_{_{J_0}}$ for some $J_0\in\esupp(S_j)$,
$j=1,2$ (see (\ref{JX})). Without loss of generality assume
$\ep_{_{J_0}}=1$ and $\ep'_{_{J_0}}=0$. Then
$\hz_{_{J_0}}^{-1}=\hu(\epb)\hz_{_{J_0}}^{-1}\hu(\epb')^{-1}\in\la
\hz_{_{\{r,s\}}},\hz_{_{J}}\mid r<s\in
J_\nu,\;J\in\J\setminus\{J_0\}\ra\ $. Set
$\b:=\theta_{j}+\tau_{_{J_0}}\in\Pi$ (see (\ref{Pi})). Therefore from
the ways $\hz_{_J}$'s and $\Pi$ are defined (see (\ref{Pi}) and
(\ref{ZJ})) and the fact that
$\hw_{\b}=\hw_{\theta_j}\hz_{_{J_0}}(\prod_{r\in
J_0}\th_{j,r})^{-1}$, it follows that $\hw_{\b}\in\la \hw_{\a}\mid
\a\in\Pi\setminus\{\b\}\ra$. But this contradicts Lemma
\ref{minimal}.\qed

\begin{rem}\label{remnew1}
{\rm By Proposition \ref{ro}, the number of integral collections for $(S_1,S_2)$ does not
depend on the particular choices of $\vd$, $S_1$ and $S_2$ in the description of $R$ in the form
(\ref{AABGP}). In fact, as isomorphic root systems have isomorphic Weyl groups, it follows from Proposition
\ref{ro} that $$\hbox{Inc}(R):=\hbox{ the number of integral collections for }
(S_1,S_2)$$ is an isomorphism invariant of
$R$. Since the trivial collection is an integral collection, we have $\hbox{Inc}(R)\geq 1$.
}
\end{rem}

\begin{cor}\label{fgag}
$\ker(\psi)$ is isomorphic to a direct sum of at most
$|\J|$-copies of $\bbbz_2$.
\end{cor}
\proof From Proposition \ref{ro}(ii) and the fact that there are
at most $2^{|\J|}$ integral collections for $(S_1, S_2)$  we have
$|\ker(\psi)|\leq2^{|\J|}$. Also by Proposition \ref{ro}(i) each
non-trivial element of $\ker(\psi)$ is of order 2, therefore the result
is clear as $\ker(\psi)$ is abelian.$\quadd$ \qed

Let $n_0\in\bbbz_{\geq 0}$ be the number of copies of $\bbbz_2$
involved in $\ker(\psi)$, then
\begin{eqnarray}\label{int-num}
|\ker(\psi)|=2^{n_0}\andd n_0=\log_2\hbox{Inc}(R).
\end{eqnarray}

\begin{cor}\label{fg}
{(i)} If $\J=\emptyset$ {(}in particular for types $X=F_4$
or $G_2${)}, then $n_0=0$.

 {(ii)} If $X=B_2$ and  $\{r,s\}\in\supp(S_1)\cup\supp(S_2)$ for all $r,s\in J_t$ or $ r,s\in J^t_\nu$, then
$n_0=|\esupp(S_1)|+|\esupp(S_2)|$. In particular, if $S_1$ and
$S_2$ are lattices, then
$$n_0=2^{\nu-t}+2^{t}-2-\nu-\left(
\begin{array}{c}
  t \\
  2 \\
\end{array}\right)-\left(
\begin{array}{c}
  \nu-t \\
  2 \\
\end{array}\right).$$

{(iii)} If $X=B_\ell$ $(\ell\geq2)$ and $\{r,s\}\in\supp(S_1)$
for all $r,s\in J_t$, then   $n_0=|\esupp(S_1)|$. In particular,
if $S_1$ is a lattice, then
$$n_0= 2^{t}-1-t-\left(
\begin{array}{c}
  t \\
  2 \\
\end{array}\right).$$

{(iv)} If $X=C_\ell$ $(\ell\geq2)$ and $\{r,s\}\in\supp(S_2)$
for all $ r,s\in J^t_\nu$, then $n_0=|\esupp(S_2)|$. In
particular, if $S_2$ is a lattice, then
$$n_0=  2^{\nu-t}-1-\nu+t-\left(
\begin{array}{c}
 \nu-t \\
  2 \\
\end{array}\right).$$
\end{cor}

\proof By \cite[Example 1.4(ii)]{AS2}, under the assumptions given in each statement, there is exactly $2^{n_0}$ integral collections. Now
the result follows from Proposition \ref{ro}(ii).\qed

\begin{DEF}\label{min-ears}
{\rm Let $R$ be a reduced extended affine root system with extended affine
Weyl group $\w$. Following \cite{H1}, we  call  $R$ a {\it
minimal} extended affine root system, if there is no
$\a\in\r^\times$ such that the reflections associated to the
elements of the set $R^\times\backslash\w\a$ generate $\w$.}
\end{DEF}

We now state our main theorem about the presentation by conjugation.
\begin{thm}\label{Theorem-3}
Let  $R=R(X, S_1, S_2)$ be a reduced non-simply laced extended
affine root system of the form { (\ref{AABGP})} with extended
affine Weyl group $\w$. Then the following statements are equivalent:

{(a)} $\w$ has the presentation by conjugation.

{(b)} $Z(\hat{\w}){\cong} Z(\w)$.

{(c)} $Z(\hat\w)$ is a free abelian group.

{(d)} The epimorphism $\psi:\hat\w\longrightarrow\w$ given by (\ref{main}) is injective.

{(e)} $\{w_{\a}\mid \a\in\Pi\}$ is a minimal set of generators
for $\w$, where $\Pi$ is given by {(\ref{Pi})}.

{(f)} $\hbox{Inc}(R)=1$.

{(g)} $R$ is a minimal   root system.
 \end{thm}
 \proof Clearly $(a)\Rightarrow$(b), and assertion (b)$\Rightarrow$(c)
 follows from Proposition \ref{induces}(iii) and (\ref{free-4}).
 From Lemma \ref{induces}(i), Proposition \ref{ro} and (\ref{free-4}) we have (c)$\Leftrightarrow$(d).
 The assertion $(d)\Rightarrow (e)$
 follows immediately from Lemma \ref{minimal}.
 We now show that $(e)\Rightarrow (d)$. Clearly   $(d)$ holds, if $X=F_4$ or
 $G_2$ (see  Corollary \ref{fg}). Next, let
$X=B_\ell$ or $C_\ell$.  Suppose to the contrary that
 (e) holds but (d) does not. Let $1\not=u\in\ker(\psi)$. Then by  Proposition
 \ref{ro}(i),
 $\hu$ can be written in the form
  $u=\hu(\ep)$, where  $\epb=\{\ep_{_J}\}_{J\in\J}$
is a nontrivial integral collection   for $R$. Let
$\ep_{_{J_0}}=1$ for some $J_0\in\J$ and set
$\b:=\theta_{j}+\tau_{_{J_0}}\in\Pi$, where $J_0\in\esupp(S_j)$,
$j=1,2$ (see (\ref{JX})). Then applying the homomorphism $\psi$ to the both sides of $\hu=\hu\big(\epb)$ we have
 $$z^{-1}_{_{J_0}}=\prod_{1\leq r<s\leq\nu}z_{_{\{r,s\}}}^{m_{r,s}}\prod_{J\in\mathcal
 J\setminus\{J_0\}} z_{_{J}}^{\epsilon_{_J}},\;\;\mbox{where $z_{_{J}}:=\psi(\hz_{_{J}}$).}$$
 Now using an argument similar to the proof of Proposition \ref{ro}(ii) we
get that $w_{\b}\in\la w_\a\mid\a\in\Pi\setminus\{\b\}\ra$ which
contradicts (e). The equivalence of (f) and (g) is immediate from
Proposition \ref{ro}. This completes the proof that the first six
assertions are equivalent. Finally the equivalence
$(d)\Leftrightarrow (g)$ is proved in \cite[Theorem 5.8]{H1}.\qed

Using Theorem \ref{Theorem-3} and Corollary \ref{fg} we have the
following result.
\begin{cor}\label{result2}
 If $\J=\emptyset$ {(}in particular for types $F_4$  and $G_2${)},
 then $\w$ has the presentation by conjugation.
\end{cor}

Let $\d_1,\d_2$ be  given by (\ref{supofu}). As in \cite[\S 1]{AS2},
we call a collection
  $\{\ep_{_J}\}_{J\in\esupp(S_j)}$,  $1\leq j\leq 2$,
  $\ep_{_J}\in\{0,1\}$,
   an {\it integral collection}  for $S_j$, if
  \begin{eqnarray*}
 \d_j(r,s)^{-1}\hspace{-.4cm}\sum_{J\in\esupp(S_j)}
 \chi_{_J}(r,s)\epsilon_{_J}\in\bbbz,
\end{eqnarray*}
for all $r<s\in J_t$ if $j=1$, and $r<s\in J_\nu^t$ if $j=2$.
Then employing the notion of integral collection for
semilattices $S_1$ and $S_2$ and using (\ref{integ-R}),
(\ref{JX}), (\ref{delta}) and Theorem \ref{Theorem-3},
 we have the following corollary.
\begin{cor}\label{FG}
(i)  If $X=B_2$,  then $\w$ has the presentation by conjugation if and only
 if  the trivial collection is the only integral collection for
$S_1$ and $S_2$.

(ii) If $X=B_\ell$${ (}\ell\geq3{ )}$, then $\w$ has the presentation by conjugation
if and only if the trivial collection
is the only integral collection for $S_1$.

(iii) If $X=C_\ell$${ (}\ell\geq3{ )}$, then $\w$ has the presentation by conjugation
if and only if the trivial
collection is the only integral collection for $S_2$.

\end{cor}

 Using Corollary \ref{FG},  \cite[Example 1.4(iii) and Lemma 1.5]{AS2} and     Theorem
\ref{Theorem-3} we have the  following corollaries:
\begin{cor}\label{Cor-B2}
If  $X=B_2$, then $R$ is a minimal extended affine root system if
one of the following conditions holds:

(a) $\ind(S_1)-t\leq 3$ and $\ind(S_2)-(\nu-t)\leq3$,

(b) $t\leq3$, $\nu-t\leq 3$, $\ind(S_1)\neq7$ and
$\ind(S_2)\neq7$.

 In particular, if
 $\nu-t\leq 3$ and $t\leq3$, then $$\mbox{$R$ is minimal $\Longleftrightarrow$ $\ind(S_1)\neq7$ and
$\ind(S_2)\neq7$.}$$
\end{cor}

\begin{cor}\label{B-ell}
If $X=B_{\ell}$ ($\ell\geq3$), then $R$ is a minimal extended
affine root system   if one of the following conditions holds:

(a) $\ind(S_1)-t\leq 3$,

(b) $t\leq3$ and $\ind(S_1)\neq7$.

In particular,  if   $t\leq3$, then $R$ is minimal if and only if
$\ind(S_1)\neq7$,
\end{cor}

\begin{cor}\label{C-ell}
 If $X=C_{\ell}$ ($\ell\geq3$), then $R$ is a minimal extended affine root
system  if one of the following conditions holds:

(a)  $\ind(S_2)-(\nu-t)\leq3$,

(b)  $\nu-t\leq3$ and $\ind(S_2)\neq7$.

In particular,  if   $\nu-t\leq3$, then $R$ is minimal if and only
if $\ind(S_2)\neq7$.
\end{cor}

\begin{rem}\label{root-index}
In \cite[(4.16)]{A4}, the author defines a notion of index, denoted $\ind(R)$, for a reduced extended affine root system $R$ and shows that if $\ind(R)=0$, then $\w$ has the presentation by conjugation. From Corollaries \ref{Cor-B2}--\ref{C-ell}, it is immediate that an extended affine root system of index zero is minimal and so by Theorem \ref{Theorem-3} has the presentation by conjugation.
\end{rem}

\begin{cor}\label{exa-B}
(i) If $X=B_\ell$ and  $S_1$ satisfies one of the
following conditions, then $R$ is not minimal:

\begin{itemize}
\item[(a)] There exists  $J\in\esupp(S_1)$ such that
$\{r,s\}\in\supp(S_1)$ for all $r,s\in J$.

\item[(b)]  $t\geq3$ and $S_1$ is a lattice.

 \item[(c)]  $t>3$ and $\ind(S_1)=2^t-2$.
\end{itemize}
 (ii)  If $X=C_\ell$ and  $S_2$ satisfies one of the following conditions then
 $R$ is not minimal:

\begin{itemize}
\item[(a)] There exists  $J\in\esupp(S_2)$ such that
$\{r,s\}\in\supp(S_2)$ for all $r,s\in J$.

\item[(b)]  $\nu-t\geq3$ and $S_2$ is a lattice.

\item[(c)] $\nu-t>3$ and $\ind(S_2)=2^{\nu-t}-2$.
\end{itemize}
\end{cor}

\begin{pro}\label{exist}
Let $X=B_\ell, C_\ell$. and $\nu, t, m_1$ and $m_2$ be integers such that $7\leq t+4\leq
m_1\leq2^t-1$ and  $7\leq\nu-t+4\leq m_2\leq2^{\nu-t}-1$. Then
there exists an extended affine root system $R=R(X,S_1,S_2)$ of nullity $\nu$,  twist number $t$ with
\begin{equation}\label{S1}
   \ind(S_1)=\left\{\begin{array}{ll}
  m_1 & \hbox{if $X=B_\ell$, }\vspace{2mm} \\
  2^t-1 & \hbox{if $X=C_{\ell\geq3}$ ,}
\end{array}\right.
\end{equation}
\andd
\begin{equation}\label{S2}\ind(S_2)=\left\{\begin{array}{ll}
  m_2 & \hbox{if $X=C_\ell$}\vspace{2mm} \\
 2^{\nu-t}-1 & \hbox{if $X=B_{\ell\geq3}$,}
\end{array}\right.
\end{equation}
such that $R$ is not minimal.
\end{pro}

\proof Let $\rd$ be a finite root system of type $B_\ell$ or
$C_\ell$ in $\vd$ equipped with a positive definite symmetric bilinear
form $\fm$. Let $\v=\vd\oplus\vz_1\oplus\vz_2$, where
$\dim\vz_1=t$ and
$\dim\vz_2=\nu-t$. We extend the form $\fm$ on $\vd$ to a
positive semi-definite symmetric bilinear form, denoted again  by
$\fm$, on $\v$ as follows
$$\fm_{_{\vd\times\vd}}:=\fm\andd(\v,\vz):=\{0\}.$$ By the proof of \cite[Corollary 5.18]{AS2}, there exist semilattices $S_1$ and $S_2$ in $\vz_1$
and $\vz_2$, respectively  satisfying (\ref{S1}) and (\ref{S2})
such that nontrivial integral collections exist for them. Then
$$R=R(X, S_1,S_2):=(S+S)\cup(\rds+S_1\oplus\la S_2\ra)\cup(\rdl+2\la S_1\ra\oplus S_2)$$
 is an extended affine root system
of type $X$, nullity $\nu$ with twist number $t$. By Theorem
\ref{Theorem-3},  $R$ is  not a minimal root system. \qed

\section{\bf Appendix}\setcounter{equation}{0}\label{Appendix}
In the proof of Proposition \ref{pox}(iii), we used
several results from literature on the concept of the presentation by conjugation. For the convenience of the reader we provide here a direct proof for part (iii) of Proposition \ref{pox} which is self-contained,
and unlike the earlier proof does not require any
previous knowledge about the existence of the presentation by conjugation for low nullities.

\begin{lem}\label{impor}
$[\hat{t}_{i,r},\hat{t}_{j,s}]=
\hz_{_{\{r,s\}}}^{\Delta(r,s)^{-1}a_{i,j}(r,s)}$, $i,j\in J_\ell$, $r\leq  s\in J_\nu$.
\end{lem}
\proof  Using  (\ref{tir}) and Lemma \ref{ip}(ii)  we have
\begin{eqnarray}\label{poxfg}
 [\th_{i,r},\;\th_{j,s}]=
\th^{-(\a_j,\a^\vee_i)k_{i,r}\sg_r}_{\a_j+k_{j,s}\sg_s}\th^{(\a_j,\a^\vee_i)k_{i,r}\sg_r}_{\a_j}.
\end{eqnarray}
If $(\a_i,\a^\vee_j)=0$, then $a_{i,j}(r,s)=0$ and so by
(\ref{poxfg}) we are done. So we may assume $(\a_i,\a_j)\not=0$.
The only cases which we must consider are:

(1) $(\a_i,\a_j)\in \rd\times\rds$ with $ r\leq s\in J_t$ or
$(\a_i,\a_j)\in  \rd\times\rdl$ with $r\leq s\in J^t_\nu$,

(2) $(\a_i,\a_j)\in \rds\times\rdl$ with $ r\leq s\in J_t$ or
$(\a_i,\a_j)\in  \rdl\times\rds$ with $r\leq s\in J^t_\nu$,

(3) $(\a_i,\a_j)\in\rds\times\rdl$ with  $(r,s)\in J_t\times
J^t_\nu$,

(4) $(\a_i,\a_j)\in\rds\times\rds$ with $ r\leq s\in J^t_\nu$ or
$(r,s)\in J_t\times J^t_\nu$,

(5) $(\a_i,\a_j)\in\rdl\times\rdl$  with $ r\leq s\in J_t$ or
    $(r,s)\in J_t\times J^t_\nu$,

 (6) $(\a_i,\a_j)\in\rdl\times\rds$ with  $(r,s)\in J_t\times
 J^t_\nu$.\\
If   (1)  holds, then by (\ref{fact-1}) and (\ref{kr}) we have
$k_{j,r}=k_{j,s}=1$ and $kk_r^{-1}\a_j^\vee=\a_j$ and so
 $a_{i,j}(r,s)=(\a_j,\a^\vee_i)k_{i,r}$ (see (\ref{aij})). First, let
$\{r,s\}\in\supp(S_1)\cup\supp(S_2)$. Then
   $\Delta(r,s)=1$ (see (\ref{delta}))
 and so using (\ref{poxfg})   we have
\begin{eqnarray*}
[\th_{i,r},\;\th_{j,s}]&=&
\th^{-a_{i,j}(r,s)\sg_r}_{\a_j+\sg_s}\th^{a_{i,j}(r,s)\sg_r}_{\a_j}\\
\mbox{(Lemmas \ref{gen-3} and
\ref{conj-cen-1}(ii))}&=&(\th^{-\sg_r}_{\a_j+\sg_s}\th^{\sg_r}_{\a_j})^{a_{i,j}(r,s)}=
((\th^{\sg_r}_{\a_j+\sg_s})^{-1}\th^{\sg_r}_{\a_j})^{a_{i,j}(r,s)}\\
\hbox{(using  Lemma \ref{gen-3})}&=&
((\th^{\sg_s+\sg_r}_{\a_j}\th^{\sg_s}_{-\a_j})^{-1}\th^{\sg_r}_{\a_j})^{a_{i,j}(r,s)}\\
\hbox{(using  Lemma \ref{gen-3}  and (\ref{ZJ}))}&=&
(\th^{\sg_s}_{\a_j}\th^{-\sg_s-\sg_r}_{\a_j}\th^{\sg_r}_{\a_j})^{a_{i,j}(r,s)}=
\hz_{_{\{r,s\}}}^{\Delta(r,s)^{-1}a_{i,j}(r,s)}.
\end{eqnarray*}
(We have used the fact that if $xy$ is central for two elements
$x$, $y$ of a group $G$, then $(xy)^n=x^ny^n$ and $xy=yx$). Next,
let $\{r,s\}\not\in\supp(S_1)\cup\supp(S_2)$. Then using
(\ref{poxfg}), (\ref{bvc}) and   the facts that
$\Delta(r,s)=(\a_j,\a^\vee_j)=2$ and $k_{j,r}=k_{j,s}=1$,
 we have
\begin{eqnarray*}
[\th_{i,r},\;\th_{j,s}]&=&
\th^{-a_{i,j}(r,s)\sg_r}_{\a_j+\sg_s}\th^{a_{i,j}(r,s)\sg_r}_{\a_j}\\
\mbox{(using Lemmas \ref{gen-3} and \ref{conj-cen-1}(ii))}&=&
(\th^{-(\a_j,\a^\vee_j)k_{j,r}\sg_r}_{\a_j+k_{j,s}\sg_s}
\th^{(\a_j,\a^\vee_j)k_{j,r}\sg_r}_{\a_j})^{\Delta(r,s)^{-1}a_{i,j}(r,s)}\\
\mbox{(using (\ref{poxfg}) and
(\ref{ZJ}))}&=&[\th_{j,r},\;\th_{j,s}]^{\Delta(r,s)^{-1}a_{i,j}(r,s)}=
   \hz_{_{\{r,s\}}}^{\Delta(r,s)^{-1}a_{i,j}(r,s)}.
\end{eqnarray*}
   If (2)
   holds,
   then    we
have $k_{i,r}=k_{i,s}=1$, $kk_r^{-1}\a_i^\vee=\a_i$ and
 $a_{i,j}(r,s)=(\a_i,\a^\vee_j)k_{j,s}$ and so using an argument similar
  to the case (1)  we get
    \begin{eqnarray*}
[\th_{j,s},\;\th_{i,r}]=
\hz_{_{\{r,s\}}}^{-\Delta(r,s)^{-1}(\a_i,\a^\vee_j)k_{j,s}\sg_s}=
\hz_{_{\{r,s\}}}^{-\Delta(r,s)^{-1}a_{i,j}(r,s)}.
\end{eqnarray*} If   (3) holds, then from the fact that $\Delta(r,s)=-a_{i,j}(r,s)=1$
and  the way that $\hz_{_{\{r,s\}}}$ is defined
 we have
\begin{eqnarray*}
[\th_{i,r},\;\th_{j,s}]=[\th_{j,s},\;\th_{i,r}]^{-1}=\hz_{_{\{r,s\}}}^{-1}=
\hz_{_{\{r,s\}}}^{\Delta(r,s)^{-1}a_{i,j}(r,s)}.
\end{eqnarray*} If
 (4) holds, then by (\ref{dotPi}) and
the assumption  $(\a_i,\a_j)\not=0$     we have
$k_{i,r}=k_{1,r}=k_{j,s}=k_{2,s}=1$ and  $(\a_i,\a^{\vee}_j)=
  (\theta_1,\a^{\vee}_2)n=-n$ where
\begin{eqnarray}\label{n}
n:=\left\{\begin{array}{ll}
  1, & \hbox{if $\a_i\neq\a_j$,}\vspace{2mm} \\
 -2, & \hbox{if $\a_i=\a_j$.}
\end{array}\right.
\end{eqnarray}  Then using (\ref{poxfg}) and the fact that $a_{i,j}(r,s)=na_{1,2}(r,s)$
    we have \begin{eqnarray*}
[\th_{i,r},\;\th_{j,s}]&=&  [\th_{j,s},\;\th_{i,r}]^{-1} =
(\th^{-(\theta_1,\theta^{\vee}_2)nk_{2,s}\sg_s}_{\a_i+k_{i,r}\sg_r}
\th^{(\theta_1,\theta^{\vee}_2)nk_{2,s}\sg_s}_{\a_i})^{-1}\\
\th^{(\theta_1,\theta^{\vee}_2)nk_{2,s}\sg_s}_{\theta_1})^{-1}\\
 \hbox{(using Lemmas \ref{well-define}(ii)) and \ref{gen-3})}&=&
 (\th^{-(\theta_1,\theta^{\vee}_2)k_{2,s}\sg_s}_{\theta_1+k_{1,r}\sg_r}
 \th^{(\theta_1,\theta^{\vee}_2)k_{2,s}\sg_s}_{\theta_1})^{-n}=
 [\th_{2,s},\;\th_{1,r}]^{-n}\\
 &=&[\th_{1,r},\;\th_{2,s}]^{n}\\
\hbox{(using case (1) or (3)) }&=&
\hz_{_{\{r,s\}}}^{\Delta(r,s)^{-1}na_{1,2}(r,s)}=
\hz_{_{\{r,s\}}}^{\Delta(r,s)^{-1}a_{i,j}(r,s)}.
\end{eqnarray*}
If
 (5) holds, then  $k_{j,s}=k_{2,s}$,\;  $(\a_j,\a^{\vee}_i)k_{i,r}=
  (\theta_2,\theta^{\vee}_1)k_{1,r}n$ where $n$ is given by (\ref{n})
    and    $a_{i,j}(r,s)=na_{1,2}(r,s)$ and so using (\ref{poxfg}) and  Lemmas \ref{well-define}(ii)
    and \ref{gen-3}
    we have
 \begin{eqnarray*}
[\th_{i,r},\;\th_{j,s}]&=&\th^{-(\theta_2,\theta^{\vee}_1)n
 k_{1,r}\sg_r}_{\theta_2+k_{2,s}\sg_s}\th^{(\theta_2,\theta^{\vee}_1)nk_{1,r}\sg_r}_{\theta_2}\\
&=&(\th^{-(\theta_2,\theta^{\vee}_1)k_{1,r}\sg_r}_{\theta_2+
k_{2,s}\sg_s}\th^{(\theta_2,\theta^{\vee}_1)k_{1,r}\sg_r}_{\theta_2})^n=
[\th_{1,r},\;\th_{2,s}]^n\\
\hbox{(using case (2) or (3)) }&=&
\hz_{_{\{r,s\}}}^{\Delta(r,s)^{-1}na_{1,2}(r,s)}=
\hz_{_{\{r,s\}}}^{\Delta(r,s)^{-1}a_{i,j}(r,s)}.
\end{eqnarray*}
If (6) holds, then  using (\ref{poxfg}) and the facts that
$(\a_j,\a^\vee_i)=-k_{j,s}=-\Delta(r,s)=-1$ and
$k_{i,r}=k_r=k=-a_{i,j}(r,s)$  we have
\begin{eqnarray*}
 [\hat{t}_{i,r}, \hat{t}_{j,s}]
&=& \th^{-\Delta(r,s)^{-1}a_{i,j}(r,s)\sg_r}_{\a_j+\sg_s}\th^{\Delta(r,s)^{-1}a_{i,j}(r,s)\sg_r}_{\a_j}\\
&=& (\th^{-k_{j,r}\sg_r}_{\a_j-(\a_j,\a^\vee_i)k_{i,s}\sg_s}
\th^{k_{j,r}\sg_r}_{\a_j})^{\Delta(r,s)^{-1}a_{i,j}(r,s)}\\
\hbox{(using  Lemma \ref{ip}(ii) and (\ref{ZJ}))} &=&
[\hat{t}_{i,s}, \hat{t}_{j,r}]^{\Delta(r,s)^{-1}a_{i,j}(r,s)}\\
\hbox{(using  case (3))}&=&
\hz_{_{\{r,s\}}}^{\Delta(r,s)^{-1}a_{i,j}(r,s)}.\hspace{3.5cm}\Box
 \end{eqnarray*}


 \end{document}